\documentclass{birkjour}
 \usepackage{mathabx}
 \usepackage{xcolor}
\numberwithin{equation}{section}

\newtheorem{theorem}{Theorem}[section]

\newtheorem{definition}[theorem]{Definition}

\newtheorem{lemma}[theorem]{Lemma}

\newtheorem{proposition}[theorem]{Proposition}

\newtheorem{rem}{Remark}[section]

\begin{document}

%
%
%
%
%
%
%
%
%

\title[Toeplitz Operators Defined between Köthe Spaces]
 {{Toeplitz Operators Defined between Köthe Spaces}}

\author{Nazl\i\;Do\u{g}an}
\address{Fatih Sultan Mehmet Vakif University, 
    34445 Istanbul, Turkey}
\email{ndogan@fsm.edu.tr}




\begin{abstract}
The aim of this paper is to define Toeplitz operators between Köthe spaces, especially power series spaces. We determine the conditions for continuity and compactness of these operators. We define the concept of S-tameness of a family of continuous operators. We construct some conditions on $S$-tameness for the families consisting of Toeplitz operators.
\end{abstract}

\maketitle
\section{Introduction}
The theory of Toeplitz operators defined on a Hilbert space such as the Hardy and  the Bergman spaces is well-studied. The significance of this theory lies in the connection between operator theory and function spaces. The matrix of a Toeplitz operator defined on Hardy space of unit disk $H^{2}(\mathbb{D})$ is a Toeplitz matrix. Moreover, any bounded operator on $H^{2}(\mathbb{D})$ whose associated matrix is a Toeplitz matrix is a Toeplitz operator. We direct the reader to \cite{Ros} for more information about Toeplitz operators defined on $H^{2}(\mathbb{D})$. In recent years, Toeplitz operators, whose "associated" matrix is Toeplitz, are defined for more general topological vector spaces. For instance, in \cite{DJ}, Doma\'nski and Jasiczak developed the analogous theory for the space of $\mathcal{A}(\mathbb{R})$ real analytic functions on the real line. This space is not a Banach space, even not a metrizable space. In \cite{J}, Jasiczak introduced and characterized the class of Toeplitz operators
on the Fr\'echet space of all entire functions $\mathcal{O}(\mathbb{C})$.

In \cite{J}, Jasiczak defined  a continuous linear operator on $\mathcal{O}(\mathbb{C})$ as a Toeplitz operator if its matrix
is a Toeplitz matrix. The matrix of an operator is defined with
respect to the Schauder basis $(z^{n})_{n\in{N_{0}}}$. In this case, the symbol of a Toeplitz operator comes from the space $\mathcal{O}(\mathbb{C})\oplus (\mathcal{O}(\mathbb{C}))^{\prime}_{b}$ here $(\mathcal{O}(\mathbb{C}))^{\prime}_{b}$ is the strong dual of $\mathcal{O}(\mathbb{C})$. The space of entire functions $\mathcal{O}(\mathbb{C})$ is isomorphic to a power series space of infinite type $\Lambda_{\infty}(n)$. By taking inspiration from Jasiczak paper \cite{J}, we will define Toeplitz operators on more general power series spaces of finite or infinite type. Similarly, we will show that the symbol space is $\Lambda_{r}(\alpha)\oplus (\Lambda_{r}(\alpha))^{\prime}$, for $r=1,\infty$.
We will also search the compactness of these operators and mention the conditions for the tameness of the families consisting of Toeplitz operators.

\section{Preliminaries}
In this section, we provide some fundamental facts and definitions essential for the subsequent discussions. We will use the standard terminology and notation of \cite{MV}.
 
A complete Hausdorff locally convex space E whose topology defined by countable fundamental
system of seminorms $(\|\cdot\|_{k})_{k\in \mathbb{N}}$ is called a Fr\'echet space. A matrix $(a_{n,k})_{k,n\in \mathbb{N}}$ of non-negative numbers is called a Köthe matrix if it is satisfies  the following conditions:
\begin{itemize}
\item[1.] For each $n\in \mathbb{N}$ there exists a $k\in \mathbb{N}$ with $a_{n,k}>0$. 
\item[2.] $a_{n,k}\leq a_{n,k+1}$ for each $n,k\in \mathbb{N}$.
\end{itemize}
For a Köthe matrix $(a_{n,k})_{n,k\in \mathbb{N}}$,
$$ K(a_{n,k})=\bigg\{ x=(x_{n})_{n\in \mathbb{N}}: \;\|x\|_{k}:=\sum^{\infty}_{n=1}|x_{n}|a_{n,k}< \infty \quad\text{for all}\quad k\in \mathbb{N} \bigg\}$$
is called a Köthe space. Every Köthe space is a Fr\'echet space. From Proposition 27.13 of \cite{MV}, the dual space of a Köthe space is
$$ (K(a_{n,k}))^{\prime}=\bigg\{ y=(y_{n})_{n\in \mathbb{N}}\; \bigg|\; \sup_{n\in \mathbb{N}}|y_{n}a_{n,k}^{-1}|<+\infty \; \text{for some}  \; k\in \mathbb{N}\bigg\}.$$

By Grothendieck-Pietsch Criteria \cite{MV}[Theorem 28.15] a Köthe space $K(a_{n,k})$ is nuclear if and only if  for every $k\in \mathbb{N}$, there exists a $l>k$ so that $$\displaystyle \sum^{\infty}_{n=1}{a_{n,k}\over a_{n,l}}<\infty.$$
For a nuclear Köthe spaces, $\|x\|_{k}=\sup_{n\in \mathbb{N}}|x_{n}|a_{n,k}$, $k\in \mathbb{N}$ forms an equivalent system of seminorms to the fundamental system of seminorms $\|x\|_{k}=\sum^{\infty}_{n=1}|x_{n}|a_{n,k}$, $k\in \mathbb{N}$, see Proposition 28.16 of \cite{MV}.
~\\

Dynin-Mitiagin Theorem \cite{MV}[Theorem 28.12] states that if a nuclear Fr\'echet
space E with the sequence of seminorms $(\|f_{n}\|)_{n\in \mathbb{N}}$ has a Schauder basis $(f_{n})_{n\in \mathbb{N}}$,
then it is canonically isomorphic to a nuclear Köthe space defined by the matrix
$(\|f_{n}\|_{k})_{n,k\in \mathbb{N}}$. So, nuclear Köthe spaces hold a significant place in theory of nuclear Fr\'echet space. 

Let $\alpha=\left(\alpha_{n}\right)_{n\in \mathbb{N}}$ be a non-negative increasing sequence with $\displaystyle \lim_{n\rightarrow \infty} \alpha_{n}=\infty$. A power series space of finite type is defined by
$$\Lambda_{1}\left(\alpha\right):=\left\{x=\left(x_{n}\right)_{n\in \mathbb{N}}: \;\left\|x\right\|_{k}:=\sum^{\infty}_{n=1}\left|x_{n}\right|e^{-{1\over k}\alpha_{n}}<\infty \textnormal{ for all } k\in \mathbb{N}\right\}$$
and a power series space of infinite type is defined by
$$\displaystyle \Lambda_{\infty}\left(\alpha\right):=\left\{x=\left(x_{n}\right)_{n\in \mathbb{N}}:\; \left\|x\right\|_{k}:=\sum^{\infty}_{n=1}\left|x_{n}\right|e^{k\alpha_{n}}<\infty \textnormal{ for all } k\in \mathbb{N}\right\}.$$
Power series spaces form an important family of Köthe spaces  and they contain the spaces of holomorphic functions on $\mathbb{C}^{d}$ and $\mathbb{D}^{d}$,
$$\mathcal{O}(\mathbb{C}^{d})\cong \Lambda_{\infty}(n^{\frac{1}{d}}) \hspace{0.6in}\text{and} \hspace{0.6in} \mathcal{O}(\mathbb{D}^{d})\cong \Lambda_{1}(n^{\frac{1}{d}})
$$
where $\mathbb{D}$ is the unit disk in $\mathbb{C}$ and $d\in \mathbb{N}$.

Let $E$ and $F$ be Fr\'echet spaces. A linear map $T:E\to F$ is called continuous 
if for every $k\in \mathbb{N}$ there exists $p\in \mathbb{N}$ 
and $C_k>0$ such that 
$$\|Tx\|_k\leq C_k\|x\|_p $$
for all $x\in E$. A linear map $T: E\to F$ is called compact
if $T(U)$ is precompact in $F$ where $U$ is a neighborhood of zero of E. 

In this paper, we fixed the symbol $e_{n}$ to denote the sequence $$(0,0,\dots,0,1,0,\dots)$$ where 1 is in the n$^{th}$ place and 0 is in the others.

We will use the following Lemma to determine the continuity and compactness of operators defined between Köthe spaces. 
\begin{lemma}\label{Crone} 
Let $K(a_{n,k})$ and $K(b_{n,k})$ be Köthe spaces.
\begin{itemize}
\item[a.] $T: K(a_{n,k}) \to K(b_{n,k})$ is a linear continuous operator if and only if for each k there exists m such that
$$\sup_{n\in \mathbb{N}} {\frac{\|Te_{n}\|_{k}}{\|e_{n}\|_{m}}<\infty }.$$
\item[b.] If $K(b_{n,k})$ is Montel, then $T: K(a_{n,k}) \to K(b_{n,k})$ is a compact operator if and only if there exists m such that for all k
$$\sup_{n\in \mathbb{N}} {\frac{\|Te_{n}\|_{k}}{\|e_{n}\|_{m}}<\infty }.$$
\end{itemize}
\end{lemma}
\begin{proof} Lemma 2.1 of \cite{CR75}.
\end{proof}

We want to note that a Fr\'echet space $E$ is Montel if each bounded set in $E$ is relatively compact and every power series space is Montel, see Theorem 27.9 of \cite{MV}.

In the next proposition, it will be shown that the continuity condition is sufficient to ensure that linear operators defined only on the basis elements are well-defined.
\begin{proposition}\label{P0} Let $K(a_{n,k})$, $K(b_{n,k})$ be Köthe spaces and  $(a_{n})_{n\in \mathbb{N}}\in K(b_{n,k})$  be a sequence. Let us define a linear map $T:K(a_{n,k})\to K(b_{n,k})$ such as
$$Te_{n}=a_{n} \hspace{0.5in}\text{and}\hspace{0.5in} Tx=\sum^{\infty}_{n=1}x_{n}Te_{n}$$
for every $n\in \mathbb{N}$ and $x=\sum^{\infty}_{n=1} x_{n}e_{n}$. If the  continuity condition 
$$\forall k\in \mathbb{N} \quad\quad \exists m\in \mathbb{N}\quad\quad\quad \sup_{n\in \mathbb{N}}\frac{ \|Te_{n}\|_{k}}{\|e_{n}\|_{m}}<\infty$$
holds, then $T$ is well-defined and continuous operator.
\end{proposition}
\begin{proof} Let  $x=\sum^{\infty}_{n=1}x_{n}e_{n}$ be an arbitrary element of  $K(a_{n,k})$ and assume that for every $k\in \mathbb{N}$ there exists a $m\in \mathbb{N}$ such that  $$\sup_{n\in \mathbb{N}} \frac{ \|Te_{n}\|_{k}}{\|e_{n}\|_{m}}<\infty.$$
Now we define the sequence $$y_{N}=\sum^{N}_{n=1} T(x_{n}e_{n})=\sum_{n=1}^{N}x_{n}Te_{n}$$ for every $N\in \mathbb{N}$. Then for every $k\in \mathbb{N}$, there exist $m\in \mathbb{N}$ and $C>0$ such that  for sufficiently large $N,M\in \mathbb{N}$, $N>M$,
$$ \|y_{N}-y_{M}\|_{k} = \|\sum^{N}_{l=M+1}x_{l}T(e_{l})\|_{k}\leq \sum^{N}_{l=M+1}|x_{l}|\|e_{l}\|_{m} \frac{\|Te_{l}\|_{k}}{\|e_{l}\|_{m}}\leq C \sum^{N}_{l=M+1}|x_{l}|\|e_{l}\|_{m}
$$
which is arbitrarily small, since $\|x\|_{m}= \sum^{\infty}_{n=1} |x_{n}|\|e_{n}\|_{m}<\infty.$ This says that the sequence $(y_{N})_{N\in \mathbb{N}}\subseteq K(b_{k,n})$ is a Cauchy sequence and therefore has a limit. Let us say the $\displaystyle \lim_{N\to \infty}y_{N}=y$. Hence, we have 
$$ y=\lim_{N\to \infty}y_{N}=\sum^{\infty}_{n=1}x_{n}T(e_{n})=T(x)\in K(b_{k,n}).$$
Therefore $T$ is well-defined from $K(a_{k,n})$ to $K(b_{k,n})$ and $T$ is continuous from Lemma \ref{Crone}. 
\end{proof}



A grading on a Fr\'echet space $E$ is a sequence of seminorms $\{\|\cdot\|_{n}:n\in \mathbb{N}\}$ which are increasing, that is,
$$\|x\|_{1}\leq \|x\|_{2}\leq \|x\|_{3}\leq \dots$$
for each $x\in E$ and which define the topology. Every Fr\'echet space admits a grading. A graded Fr\'echet space is Fr\'echet space with a choice of grading. For more information, see \cite{H}. In this paper, unless stated otherwise, we will assume that all Fr\'echet spaces are graded Fr\'echet spaces. 

A pair of graded Fr\'echet spaces $(E,F)$ is called tame provided that there exists an increasing function $\sigma:\mathbb{N}\to \mathbb{N}$, such that for every continuous linear operator T from $E$ to $F$, there exists an $N$ and $C>0$ such that
$$\hspace{1.15in} \|Tx\|_{n} \leq C\|x\|_{\sigma(n)} \hspace{1.05in}\forall x\in E\;\; \text{and}\;\; n\geq N. $$ A Fr\'echet space $E$ is called tame if $(E,E)$ is tame. Tameness gives us a kind of control on the continuity of the operators. In \cite{DV}, Dubinsky and Vogt used this notion to find a basis in complemented subspaces of some infinite type power series spaces. We refer to \cite{A} and \cite{P} for the other studies about tameness. 

In this paper, instead of all operators defined on a Fr\'echet space, we will restrict ourselves to a subfamily of operators. With this aim, we will give the following definition.

\begin{definition} Let  $S:\mathbb{N}\to \mathbb{N}$ be a non-decreasing function. A family of linear continuous operators $\mathcal{A}\subseteq L(E,F)$ is called S-tame if for every operator $T\in \mathcal{A}$ there exist $k_{0}\in \mathbb{N}$ and  $C>0$ such that
$$\hspace{1.65in}
\|Tx\|_{k}\leq C\|x\|_{S(k)} \hspace{1in}\forall x\in E, k\geq k_{0}. $$
\end{definition}
S-tameness of a family can be given by considering only elements of bases similar to Lemma \ref{Crone}.

\begin{lemma}
Let $K(a_{n,k})$ and $K(b_{n,k})$ be K\"{o}the spaces and $S:\mathbb{N}\to \mathbb{N}$ be a non-decreasing function. A family of linear continuous operators $\mathcal{A}\subseteq L(K(a_{n,k}),K(b_{n,k}))$ is S-tame if and only if for every operator $T\in \mathcal{A}$ there exist $k_{0}\in \mathbb{N}$ and $C>0$ such that
$$\hspace{1.4in}\|Te_{n}\|_{k}\leq C
\|e_{n}\|_{S(k)} \hspace{1.1in}\forall n\in \mathbb{N}, k\geq k_{0}.$$
\end{lemma}
\begin{proof} Let us assume that for every $T\in \mathcal{A}$ there exist a $k_{0}\in \mathbb{N}$ and $C>0$ such that
$$\hspace{1.4in} \|Te_{n}\|_{k}\leq C
\|e_{n}\|_{S(k)} \hspace{1.1in}\forall n\in \mathbb{N},k\geq k_{0}.$$
Then for every $x=\sum^{\infty}_{n=1} x_{n}e_{n}$, we can write
\begin{equation*}
 \begin{split}
\|Tx\|_{k}&=\bigg\|\sum^{\infty}_{n=1} x_{n}T(e_{n})\bigg\|_{k}\leq \sum^{\infty}_{n=1}|x_{n}|\|Te_{n}\|_{k} \leq C\sum^{\infty}_{n=1} |x_{n}|\|e_{n}\|_{S(k)}\\ &=C\sum^{\infty}_{n=1} |x_{n}|a_{n,S(k)} =C\|x\|_{S(k)}.
\end{split}
\end{equation*}
This says that $\mathcal{A}$ is S-tame. The other direction is straightforward.
\end{proof}

In this paper, we will call an operator which is defined between Köthe spaces as a Toeplitz operator if its matrix is a Toeplitz matrix defined with respect to the Schauder basis $(e_{n})_{n\in \mathbb{N}}$. In sections 3 and 4, we will concentrate on lower triangular Toeplitz matrices and upper triangular Toeplitz matrices, respectively. We will determine the continuity and compactness of the operators whose associated matrix is related to these matrices. By collecting the results obtained in sections 3 and 4, we share the results about the continuity and compactness of a Toeplitz operator defined between power series spaces in section 5.
\section{Operators Defined by A Lower Triangular Toeplitz Matrix}

Let $\theta=(\theta_{n})_{n\in \mathbb{N}}$ be any sequence. The lower triangular Toeplitz matrix defined by $\theta$ is
$$
\begin{pmatrix}
\theta_{0} &0&0&0&\cdots \\
\theta_{1}&\theta_{0}&0&0&\cdots\\
\theta_{2}&\theta_{1}&\theta_{0}&0&\cdots\\
\theta_{3}&\theta_{2}&\theta_{1}&\theta_{0}&\cdots\\
\vdots&\vdots&\vdots&\vdots&\ddots
\end{pmatrix}.
$$
We want to define an operator $\widehat{T}_{\theta}: K(a_{n,k})\to K(b_{n,k})$ by taking  $\widehat{T}_{\theta}e_{n}$ as the n$^{th}$ column of the above matrix, that is,
$$\widehat{T}_{\theta}e_{n}=(0,\cdots,0,\theta_{0},\theta_{1},\theta_{2},\cdots)=\sum^{\infty}_{j=n}\theta_{j-n}e_{j}$$
provided that $\widehat{T}_{\theta}e_{n}\in K(b_{n,k})$ for every $n\in \mathbb{N}$. Therefore, for every $x=\sum^{\infty}_{n=1}x_{n}e_{n}\in K(a_{n,k})$, the operator  $\widehat{T}_{\theta}$ is given by 
$$ \widehat{T}_{\theta}x=\sum^{\infty}_{n=1}x_{n} \widehat{T}_{\theta}e_{n}.$$

In fact, we can not guarantee that the operator $\widehat{T}_{\theta}$ is correctly defined between two general K\"{o}the spaces $K(a_{n,k})$ and $K(b_{n,k})$, because we do not know if the series  $\sum^{\infty}_{n=1}x_{n} \widehat{T}_{\theta}e_{n}$ converges in $K(b_{n,k})$ for every $x\in K(a_{n,k})$
in this general case. Below, we will examine the cases where this operator can be correctly defined, and in these cases, we will analyze the continuity and compactness of this operator. Throughout this paper, we will assume that the
sequences $\theta$ satisfy the following condition
$$\forall n\in \mathbb{N} \hspace{0.2in}\exists s>n \hspace{0.6in} \theta_{s}\neq 0.$$
Otherwise, $\theta$ produces a finite rank operator, which is obviously continuous and compact. 

The idea of Proposition \ref{P1} is the same as \cite[Theorem 2.1]{UY19}. Therein one can also find a generalized form of Theorem \ref{P2}.   
\begin{proposition}\label{P1}
Let $K(a_{n,k})$, $K(b_{n,k})$ be Köthe spaces and assume that $\theta$ be a sequence and $s_{0}=\min\{t:\theta_{t}\neq 0\}$. If $\widehat{T}_{\theta}: K(a_{n,k})\to K(b_{n,k})$ is a continuous operator, then $\theta\in K(b_{n,k})$ and the following holds 
$$\hspace{0.6in}\forall k\in \mathbb{N}\quad \exists m\in \mathbb{N}, C>0 \quad \quad \quad \quad b_{n+s_{0},k}\leq Ca_{n,m} \hspace{.5in} \forall n\in \mathbb{N}.$$
\end{proposition}
\begin{proof} Let $\widehat{T}_{\theta} :K(a_{n,k}) \to K(b_{n,k})$ be a continuous operator.  For every $n\in \mathbb{N}$,
$$\widehat{T}_{\theta}e_{n}=(0,\dots,0,\theta_{0},\theta_{1},\theta_{2},\dots)=\sum^{\infty}_{j=n}\theta_{j-n}e_{j}\in K(b_{n,k}).$$
Since $\widehat{T}_{\theta} e_{1}\in K(b_{n,k})$, this gives us $\theta\in K(b_{n,k})$. By Lemma \ref{Crone}, for all $k\in \mathbb{N}$ there exist $m\in \mathbb{N}$ and $C_{1}>0$ such that
\[\hspace{.7in}\|\widehat{T}_{\theta}e_{n}\|_{k}=\sum^{\infty}_{j=n}|\theta_{j-n}|b_{j,k}\leq C_{1}\|e_{n}\|_{m}=C_{1}a_{n,m}\hspace{.7in}\forall n\in \mathbb{N}.\]
Then, for all n and $j\geq n$, we have $|\theta_{j-n}|b_{j,k}\leq C_{1}a_{n,m}.$ Hence we can write that 
$b_{n+s_{0},k} \leq C_{2}a_{n,m}$ for some $C_{2}>0$.
This completes the proof.
\end{proof}

For compact $\widehat{T}_{\theta}:K(a_{n,k})\to K(b_{n,k})$ operators, the relationship between K\"{o}the matrices $(a_{n,k})_{n,k\in \mathbb{N}}$ and $(b_{n,k})_{n,k\in \mathbb{N}}$
is as follows:
\begin{proposition}\label{P1C}
Let $K(a_{n,k})$ and $K(b_{n,k})$ be Köthe spaces such that $K(b_{n,k})$ is Montel and assume that $\theta\in K(b_{n,k})$ and $s_{0}=\min\{t:\theta_{t}\neq 0\}$. If $\widehat{T}_{\theta}: K(a_{n,k})\to K(b_{n,k})$ is a compact operator, then the following holds 
$$\hspace{.5in}\exists m\in \mathbb{N}\quad \forall k\in \mathbb{N}\quad \exists C>0  \quad \quad \quad \quad b_{n+s_{0},k}\leq Ca_{n,m} \hspace{.5in} \forall n\in \mathbb{N}.$$
\end{proposition}
\begin{proof}  The proof is similar to the proof of Proposition \ref{P1}.
\end{proof}

The converse of Proposition \ref{P1} is true when $K(b_{n,k})$ is a power series space of finite type.
\begin{theorem}\label{P2}
Let $K(a_{n,k})$ be a Köthe space, $\Lambda_{1}(\beta)$ be a power series space of finite type and assume that $\theta$ be a sequence and $s_{0}=\min\{t:\theta_{t}\neq 0\}$. $\widehat{T}_{\theta}:K(a_{n,k})\to \Lambda_{1}(\beta)$  is well-defined and continuous if and only if $\theta\in \Lambda_{1}(\beta)$ and the following condition holds:
\begin{equation}\label{Eq1}
\hspace{.2in}\forall k\in \mathbb{N} \quad \exists m\in \mathbb{N}, C>0\quad \quad \quad e^{-\frac{1}{k}\beta_{n+s_{0}}}\leq C a_{n,m}\hspace{.5in} \forall n\in \mathbb{N}.
\end{equation}
\end{theorem}
\begin{proof} If the operator $\widehat{T}_{\theta}$ is continuous, Proposition \ref{P1} gives us that $\theta\in \Lambda_{1}(\beta)$ and \eqref{Eq1} holds. For the other direction, let us assume $\theta\in \Lambda_{1}(\beta)$. By using the condition (\ref{Eq1}) we have the following: for every $k\in \mathbb{N}$, there exist $m\in \mathbb{N}$ and $C_{1}, C_{2}>0$ such that
\begin{align*}
\begin{split}
\|\widehat{T}_{\theta}e_{n}\|_{k} &\leq \|\widehat{T}_{\theta}e_{n}\|_{2k}=\sum^{\infty}_{j=n} |\theta_{j-n}|e^{-\frac{1}{2k}\beta_{j}} 
\\ &=\sum^{\infty}_{j=n+s_{0}} |\theta_{j-n}|e^{-\frac{1}{2k}\beta_{j}}e^{\frac{1}{k^{2}+2k}\beta_{n+s_{0}}}e^{-\frac{1}{k^{2}+2k}\beta_{n+s_{0}}}\\& 
 \leq C_{1}a_{n,m} \sum^{\infty}_{j=n+s_{0}} |\theta_{j-n}|e^{-\frac{1}{2k}\beta_{j}}e^{\frac{1}{k^{2}+2k}\beta_{n+s_{0}}} \\ &= C_{1}a_{n,m} \sum^{\infty}_{j=n+s_{0}} |\theta_{j-n}|e^{-\frac{1}{2(k+2)}\beta_{j}}e^{-\frac{1}{k^{2}+2k}(\beta_{j}-\beta_{n+s_{0}})} \\&
 \leq C_{1} a_{n,m}\sum^{\infty}_{j=n+s_{0}} |\theta_{j-n}|e^{-\frac{1}{2(k+2)}\beta_{j}} =C_{1}\|\theta\|_{2(k+2)}a_{n,m}=C_{2}a_{n,m}
\end{split}
\end{align*}
for every $n\in \mathbb{N}$.
In the second line of the inequalites we use the fact that
$$\frac{1}{k^{2}+2k}=\frac{1}{2k}-\frac{1}{2(k+2)}$$
for every $k\in \mathbb{N}$.
Therefore, $\widehat{T}_{\theta}e_{n}\in \Lambda_{1}(\beta)$ for every $n\in \mathbb{N}$ and for every $k\in \mathbb{N}$ there exist $m\in \mathbb{N}$ such that
$$ \sup_{n\in \mathbb{N}}\frac{\|\widehat{T}_{\theta}e_{n}\|_{k}}{\|e_{n}\|_{m}}< \infty,$$
that is, $\widehat{T}_{\theta}:K(a_{n,k})\to \Lambda_{1}(\beta)$  is well defined and continuous by Proposition \ref{P0}.
\end{proof}

We can characterize the compactness of the operators $\widehat{T}_{\theta}: K(a_{n,k})\to \Lambda_{1}(\beta)$ as the follows.  
\begin{theorem}\label{P3} Let $K(a_{n,k})$ be a Köthe space, $\Lambda_{1}(\beta)$ be a power series space of finite type and assume that $\theta$ be a sequence and $s_{0}=\min\{t:\theta_{t}\neq 0\}$. $\widehat{T}_{\theta}:K(a_{n,k})\to \Lambda_{1}(\beta)$  is compact if and only if $\theta\in \Lambda_{1}(\beta)$ and the following condition holds:
\begin{equation}\label{Eq2}
\hspace{.2in}\exists m\in \mathbb{N}  \quad \forall k\in \mathbb{N} \quad \exists C>0 \quad\quad\quad e^{-\frac{1}{k}\beta_{n+s_{0}}}\leq C a_{n,m}\hspace{.3in} \forall n\in \mathbb{N}.
\end{equation}
\end{theorem}
\begin{proof} If $\widehat{T}_{\theta}$ is compact, then $\widehat{T}_{\theta}$ is continuous and $\theta\in \Lambda_{1}(\beta)$ from Proposition \ref{P1}. Proposition \ref{P1C} says that the condition \eqref{Eq2} holds.
For the other direction, assume that $\theta\in \Lambda_{1}(\beta)$ and the condition (\ref{Eq2}) holds. By  the same calculation of the proof of Theorem \ref{P2}, we can write that there exists an $m\in \mathbb{N}$ such that for all $k\in \mathbb{N}$
$$ \sup_{n\in \mathbb{N}}\frac{\|\widehat{T}_{\theta}e_{n}\|_{k}}{\|e_{n}\|_{m}}< \infty.$$
From Lemma \ref{Crone}, we can say that $\widehat{T}_{\theta}:K(a_{n,k})\to \Lambda_{1}(\beta)$ is compact.
\end{proof}

\begin{proposition}\label{P4}
Let $K(a_{n,k})$ be a Köthe space and $\Lambda_{1}(\beta)$ be a power series space of finite type. Assume there exists a non-decreasing function $S: \mathbb{N}\to \mathbb{N}$ such that the following condition holds:
\begin{equation}\label{Eq1.3}
\hspace{.2in}\forall k\in \mathbb{N} \quad \exists C>0\quad \quad \quad\quad e^{-\frac{1}{k}\beta_{n}}\leq C a_{n,S(k)}\hspace{.8in} \forall n\in \mathbb{N}.
\end{equation}
Then, the family of the operators $\widehat{T}_{\theta}:K(a_{n,k})\to \Lambda_{1}(\beta)$  for $\theta\in \Lambda_{1}(\beta)$, that is,
$$\mathcal{\widehat{A}}_{1}= \{\widehat{T}_{\theta}: \theta\in \Lambda_{1}(\beta)\}$$
is $\tilde{S}$-tame, where $\tilde{S}(k)=S(2k)$ for every $k\in \mathbb{N}$.
\end{proposition}
\begin{proof} Let $\theta\in \Lambda_{1}(\beta)$. By using \eqref{Eq1.3}, for every $k\in \mathbb{N}$, there exist $C_{1}, C_{2}>0$ such that
\begin{align*}
\begin{split}
\|\widehat{T}_{\theta}e_{n}\|_{k}&=\sum^{\infty}_{j=n} |\theta_{j-n}|e^{-\frac{2}{2k}\beta_{j}}=\sum^{\infty}_{j=n} |\theta_{j-n}|e^{-\frac{1}{2k}\beta_{j}}e^{-\frac{1}{2k}\beta_{j}} \\ &\leq \bigg(\sum^{\infty}_{j=n} |\theta_{j-n}|e^{-\frac{1}{2k}\beta_{j}}\bigg)e^{-\frac{1}{2k}\beta_{n}} \leq C_{1} \bigg(\sum^{\infty}_{j=n} |\theta_{j-n}|e^{-\frac{1}{2k}\beta_{j}}\bigg) a_{n,S(2k)} \\
& \leq C_{1} \bigg(\sum^{\infty}_{j=n} |\theta_{j-n}|e^{-\frac{1}{2k}\beta_{j-n+1}}\bigg) \|e_{n}\|_{S(2k)} \leq C_{1} \|\theta\|_{2k}\|e_{n}\|_{\tilde{S}(k)} \leq C_{2}\|e_{n}\|_{\tilde{S}(k)}.
\end{split}
\end{align*}
Therefore, for every $\theta\in \Lambda_{1}(\beta)$  and for all $k\in \mathbb{N}$, there exists $C>0$ so that we write
$$ \|\widehat{T}_{\theta}e_{n}\|_{k}\leq C \|e_{n}\|_{\tilde{S}(k)}, $$
that is, the family $\mathcal{\widehat{A}}_{1}$ is $\tilde{S}$-tame.
\end{proof}

Now we turn our attention to the power series space of infinite type. For the converse of Proposition \ref{P1}, we need the stability condition on the sequence $\beta$ when $K(b_{n,k})$ is power series space of infinite type $\Lambda_{\infty}(\beta)$. A sequence $\beta$ is called stable if
$$
\sup_{n\in \mathbb{N}} {\beta_{2n}\over \beta_{n}}<\infty.$$

\begin{theorem}\label{P5}
Let $\beta=\left(\beta_{n}\right)_{n\in \mathbb{N}}$ be a  stable sequence. $\widehat{T}_{\theta}:K(a_{n,k})\to \Lambda_{\infty}(\beta)$ is well-defined and continuous if and only if  $\theta\in \Lambda_{\infty}(\beta)$ and the following condition holds:
\begin{equation}\label{Eq3}
\hspace{.2in}\forall k\in \mathbb{N} \quad \exists m\in \mathbb{N}, C>0 \quad\quad\quad\quad e^{k\beta_{n}}\leq C a_{n,m}\hspace{.7in}\forall n\in \mathbb{N}.
\end{equation}
\end{theorem}
\begin{proof} Let assume $\theta$ be a sequence and $s_{0}=\min\{t:\theta_{t}\neq 0\}.$ If the operator $\widehat{T}_{\theta}$ is continuous, Proposition \ref{P1} gives us that $\theta\in \Lambda_{\infty}(\beta)$ and for every $k\in \mathbb{N}$, there exist $m\in \mathbb{N}$, $C>0$ such that
$$e^{k\beta_{n+s_{0}}}\leq Ca_{n,m}$$
and then we have
$$e^{k\beta_{n}}\leq Ca_{n,m}$$
for every $n\in \mathbb{N}$. This says that 
\eqref{Eq3} holds since $\beta$ is an increasing sequence. For the converse, we assume that $\beta$ is stable. Then there exists an $M> 1$ such that
$$\hspace{1.65in}\beta_{2n}\leq M\beta_{n}\hspace{1.75in} \forall n\in \mathbb{N}.$$
Since $\beta$ is increasing we have the following: if $j=2t\geq 2n$, 
$$ \beta_{j}=\beta_{2t}\leq M\beta_{t}\leq M\beta_{2t-n}=M\beta_{j-n}\leq M\beta_{j-n+1}$$
and if $j=2t+1\geq 2n$, we have $t+1\geq n$ and 
$$\beta_{j}=\beta_{2t+1}\leq \beta_{2t+2}\leq M\beta_{t+1}\leq M\beta_{2t+2-n}=M\beta_{j-n+1}.$$
Therefore, we have 
$$\beta_{j}\leq M\beta_{j-n+1}$$
for all $j\geq 2n$. Now let us assume that $\theta\in \Lambda_{\infty}(\beta)$ and the condition (\ref{Eq3}) holds. Then, for every $k\in \mathbb{N}$ there exist $m\in \mathbb{N}$, $C_{1},C_{2}>0$ such that
\begin{align*}
\begin{split}
\|\widehat{T}_{\theta}e_{n}\|_{k}&=\sum^{2n-1}_{j=n} |\theta_{j-n}|e^{k\beta_{j}}+\sum^{\infty}_{j=2n} |\theta_{j-n}|e^{k\beta_{j}} \\ &\leq
\sum^{2n-1}_{j=n} |\theta_{j-n}|e^{k\beta_{2n}} +\sum^{\infty}_{j=2n} |\theta_{j-n}|e^{Mk\beta_{j-n+1}}e^{k\beta_{j}-Mk\beta_{j-n+1}}\\
&\leq e^{Mk\beta_{n}}
\sum^{2n-1}_{j=n} |\theta_{j-n}|+\sum^{\infty}_{j=2n} |\theta_{j-n}|e^{Mk\beta_{j-n+1}} \leq e^{Mk\beta_{n}}\sum^{\infty}_{j=n} |\theta_{j-n}|e^{Mk\beta_{j-n+1}}\\
&\leq C_{1} a_{n,m}\sum^{\infty}_{j=n} |\theta_{j-n}|e^{Mk\beta_{j-n+1}}\leq C_{1}\|\theta\|_{Mk} \|e_{n}\|_{m} 
\leq C_{2}\|e_{n}\|_{m}.
\end{split}
\end{align*}
holds for every $n\in \mathbb{N}$. We employed that $\beta_{j}\leq M\beta_{j-n+1}$ for every $j\geq 2n$ in the first line. 
Therefore, $\widehat{T}_{\theta}e_{n}\in \Lambda_{\infty}(\beta)$ for every $n\in \mathbb{N}$ and for all $k\in \mathbb{N}$ there exist $m\in \mathbb{N}$ such that
$$\label{T1}
\sup_{n\in \mathbb{N}}\frac{\|\widehat{T}_{\theta}e_{n}\|_{k}}{\|e_{n}\|_{m}}< \infty
$$
that is, $\widehat{T}_{\theta}$ is well-defined and continuous by Proposition \ref{P0}.
\end{proof}
We can characterize the compactness of the operators $\widehat{T}_{\theta}: K(a_{n,k})\to \Lambda_{\infty}(\beta)$ as the follows. 
\begin{theorem}\label{P6}
Let $\beta=\left(\beta_{n}\right)_{n\in \mathbb{N}}$ be a stable sequence. $\widehat{T}_{\theta}:K(a_{n,k})\to \Lambda_{\infty}(\beta)$  is compact if and only if $\theta\in \Lambda_{\infty}(\beta)$ and the following condition holds:
\begin{equation}\label{Eq4}
\hspace{.3in}\exists m\in \mathbb{N} \quad \forall k\in \mathbb{N}\quad \exists  C>0\quad \quad\quad\quad e^{k\beta_{n}}\leq C a_{n,m}\hspace{.3in}\forall n\geq k.
\end{equation}
\end{theorem}
\begin{proof}  Let assume $\theta$ be a sequence and $s_{0}=\min\{t:\theta_{t}\neq 0\}.$ 
 If $\widehat{T}_{\theta}$ is compact, then $\theta\in \Lambda_{1}(\beta)$ from Proposition \ref{P1}. Proposition \ref{P1C}  says that there exists a $m\in N$ such that for all $k\in \mathbb{N}$ there exists a $C>0$ such that
$$e^{ k\beta_{n+s_{0}}}\leq Ca_{n,m}$$
for every $n\in \mathbb{N}$. This says that the condition (\ref{Eq4}) is satisfied since $\beta$ is an increasing sequence. 

For the converse, assume that $\theta\in \Lambda_{\infty}(\beta)$ and (\ref{Eq4}) holds. By  the same calculation of the proof of Theorem \ref{P5}, we can write that there exists $m\in \mathbb{N}$ such that for all $k\in \mathbb{N}$
$$ \sup_{n\geq k}\frac{\|\widehat{T}_{\theta}e_{n}\|_{k}}{\|e_{n}\|_{m}}< \infty$$
and then 
$$ \sup_{n\in \mathbb{N}}\frac{\|\widehat{T}_{\theta}e_{n}\|_{k}}{\|e_{n}\|_{m}}< \infty.$$
From Lemma \ref{Crone}, we can say that $\widehat{T}_{\theta}:K(a_{n,k})\to \Lambda_{\infty}(\beta)$ is compact.
\end{proof}

\begin{proposition}\label{P7}
Let $\beta=\left(\beta_{n}\right)_{n\in \mathbb{N}}$ be a non-negative increasing stable sequence with $\displaystyle \lim_{n\rightarrow \infty} \beta_{n}=\infty$. Assume that there exists a non-decreasing function $S:\mathbb{N}\to \mathbb{N}$ such that
\begin{equation}\label{Eq5}
\hspace{.2in}\forall k\in \mathbb{N}\quad \exists  C>0\quad \quad \quad\quad\quad e^{k\beta_{n}}\leq C a_{n,S(k)}\hspace{.8in} \forall n\in \mathbb{N}.
\end{equation}
Then, the family of operators $\widehat{T}_{\theta}:K(a_{n,k})\to \Lambda_{\infty}(\beta)$ for $\theta\in \Lambda_{\infty}(\beta)$
$$\mathcal{\widehat{A}}_{\infty}=\{T_{\theta}: \theta\in \Lambda_{\infty}(\beta)\}$$
is $\tilde{S}$-tame, where $\tilde{S}(k)=S(Mk)$ for every $k\in \mathbb{N}$ and here $\displaystyle M=\sup_{n\in \mathbb{N}}\frac{\alpha_{2n}}{\alpha_{n}}$.
\end{proposition}
\begin{proof} The proof is similar to the proof of Theorem \ref{P6}. Let $\theta\in \Lambda_{\infty}(\beta)$. In the proof of Theorem \ref{P5} we showed that for all $k,n\in \mathbb{N}$  
$$
\|\widehat{T}_{\theta}e_{n}\|_{k}
\leq e^{Mk\beta_{n}}\sum^{\infty}_{j=n} |\theta_{j-n}|e^{Mk\beta_{j-n+1}} \leq e^{Mk\beta_{n}}\|\theta\|_{Mk}.$$
Using the condition (\ref{Eq5}), we can write that for all $k\in \mathbb{N}$ there exists a $C>0$ such that
$$
\|\widehat{T}_{\theta}e_{n}\|_{k}
\leq e^{Mk\beta_{n}}\|\theta\|_{Mk}\leq Ca_{n, S(Mk)}=Ca_{n,\tilde{S}(k)}.$$
Then we have the following
\begin{equation*}\label{}
\hspace{1in}\forall k \quad \exists C>0 \hspace{0.7in} \|\widehat{T}_{\theta}e_{n}\|_{k}\leq C\|e_{n}\|_{\tilde{S}(k)} \hspace{1in}
\end{equation*}
for every $\theta\in \Lambda_{\infty}(\beta)$. This says that the family $\mathcal{\widehat{A}}_{\infty}$ is $\tilde{S}$-tame.
\end{proof}

\section{Operators Defined by an Upper Triangular Toeplitz Matrix}

Let $\theta=(\theta_{n})_{n\in \mathbb{N}}$ be any sequence. The upper triangular Toeplitz matrix defined by $\theta$ is
$$
\begin{pmatrix}
\theta_{0} &\theta_{1}&\theta_{2}&\theta_{3}&\cdots \\
0&\theta_{0}&\theta_{1}&\theta_{2}&\cdots\\ 0&0&\theta_{0}&\theta_{1}&\cdots\\
0&0&0&\theta_{0}&\cdots\\
\vdots&\vdots&\vdots&\vdots&\ddots
\end{pmatrix}
$$
We want to define a linear operator $\widecheck{T}_{\theta}: K(a_{n,k})\to K(b_{n,k})$ by taking  $\widecheck{T}_{\theta}e_{n}$ as the n$^{th}$ column of the above matrix, that is,
$$\widecheck{T}_{\theta}e_{n}=(\theta_{n-1},\theta_{n-2},\dots, \theta_{1},\theta_{0},0,0,\dots)=\sum^{n}_{j=1}\theta_{n-j}e_{j}.$$
Therefore, for every $x=\sum^{\infty}_{n=1}x_{n}e_{n}\in K(a_{n,k})$, the operator  $\widecheck{T}_{\theta}$ is given by 
$$ \widecheck{T}_{\theta}x=\sum^{\infty}_{n=1}x_{n} \widecheck{T}_{\theta}e_{n}.$$

In fact, we can not guarantee that the operator $\widecheck{T}_{\theta}$ is correctly defined between two general K\"{o}the spaces $K(a_{n,k})$ and $K(b_{n,k})$, because we do not know if the series  $\sum^{\infty}_{n=1}x_{n} \widecheck{T}_{\theta}e_{n}$ converges in $K(b_{n,k})$ for every $x\in K(a_{n,k})$
in this general case. Below, we will examine the cases where this operator can be correctly defined, and in these cases, we will analyze the continuity and compactness of this operator. Similiar to the previous section, we will assume that the
sequences $\theta$ satisfy the following condition
$$\forall n\in \mathbb{N} \hspace{0.2in}\exists s>n \hspace{0.4in} \theta_{s}\neq 0.$$
Otherwise, $\theta$ produces a finite rank operator, which is obviously continuous and compact. 

In \cite[Theorem 2.2]{UY19}, you can find the same idea for Proposition \ref{P8} and a generalized form of  Theorem \ref{P9} for $G_{\infty}$-spaces.
\begin{proposition}\label{P8}
Let $K(a_{n,k})$, $K(b_{n,k})$ be Köthe spaces and assume that $\theta$ be a sequence and $s_{0}=\min\{t:\theta_{t}\neq 0\}$. If $\widecheck{T}_{\theta} : K(a_{n,k})\to K(b_{n,k})$ is continuous, then $\theta\in (K(a_{n,k}))^{\prime}$ and the following holds:
$$\hspace{.25in}\forall k\in \mathbb{N}\quad \exists m\in \mathbb{N}, C>0 \quad \quad\quad\quad\quad b_{n-s_{0},k}\leq Ca_{n,m} \hspace{.65in}\forall n>s_{0}.$$
\end{proposition}
\begin{proof} Let $\widecheck{T}_{\theta} :K(a_{n,k}) \to K(b_{n,k})$ be a continuous operator. Then by Lemma \ref{Crone}, for all $k\in \mathbb{N}$ there exist $m\in \mathbb{N}$ and $C_{1}>0$ such that
$$\|\widecheck{T}_{\theta} e_{n}\|_{k}=\sum^{n}_{j=1}|\theta_{n-j}|b_{j,k}\leq C_{1}a_{n,m}\quad \quad \forall n\in \mathbb{N}.$$
Then, for all n and $j\leq n$
$$|\theta_{n-j}|b_{j,k}\leq C_{1}a_{n,m}.$$
Since $(b_{n,k})$ is a Köthe matrix, there exists a  $k_{0}\in\mathbb{N}$ satisfying $b_{1,k_{0}}\neq 0$. Let us take $C_{2}={C_{1}\over b_{1,k_{0}} } $. Then we can write 
$$|\theta_{n-1}|\leq C_{2} a_{n,m}$$
and this says that $\theta\in (K(a_{n,k}))^{\prime}$.
Further, we have that for every $k\in \mathbb{N}$, there exist $m\in \mathbb{N}$ and $C
_{1}>0$ so that 
$$|\theta_{s_{0}}|b_{n-s_{0}}\leq C_{1}a_{n,m}$$
and 
$$b_{n-s_{0}}\leq \frac{C_{1}}{\theta_{s_{0}}} a_{n,m}.$$
hold for all $n> s_{0}$.
By choosing $\displaystyle C=\frac{C_{1}}{|\theta_{s_{0}}|}$, we have $$\hspace{1.5in}b_{n-s_{0},k}\leq C_{3}a_{n,m} \hspace{1.5in}\forall n>s_{0}.$$
This completes the proof.
\end{proof}

For compact $\widecheck{T}_{\theta}:K(a_{n,k})\to K(b_{n,k})$ operators, the relationship between K\"{o}the matrices $(a_{n,k})_{n,k\in \mathbb{N}}$ and $(b_{n,k})_{n,k\in \mathbb{N}}$
is as follows:
\begin{proposition}\label{P8C}
Let $K(a_{n,k})$ and $K(b_{n,k})$ be Köthe spaces such that $K(b_{n,k})$ is Montel and assume that $\theta\in K(b_{n,k})$ and $s_{0}=\min\{t:\theta_{t}\neq 0\}$. If $\widecheck{T}_{\theta}: K(a_{n,k})\to K(b_{n,k})$ is a compact operator, then the following holds 
$$\hspace{.35in}\exists m\in \mathbb{N}\quad \forall k\in \mathbb{N}\quad \exists C>0 \quad \quad \quad \quad b_{n-s_{0},k}\leq Ca_{n,m} \hspace{.5in} \forall n>s_{0}.$$
\end{proposition}
\begin{proof}  The proof is similar to the proof of Proposition \ref{P8}.
\end{proof}

The converse of Proposition \ref{P1} is true when $K(b_{n,k})$ is a power series space of finite type.

\begin{theorem}\label{P9} Let $\Lambda_{\infty}(\alpha)$ be a nuclear power series space of infinite type, $K(b_{n,k})$ be a Köthe space and assume that $\theta$ be a sequence and $s_{0}=\min\{t:\theta_{t}\neq 0\}$. 
The operator $\widecheck{T}_{\theta}:\Lambda_{\infty}(\alpha)\to K(b_{n,k})$ is well-defined and continuous if and only if $\theta\in (\Lambda_{\infty}(\alpha))^{\prime}$  and the following condition holds:
\begin{equation}\label{Eq6}
\hspace{0.2in} \forall k\in \mathbb{N} \quad \exists m\in \mathbb{N}, C>0\quad \quad\quad \quad b_{n-s_{0},k}\leq C e^{m\alpha_{n}} \hspace{0.4in}\forall n>s_{0}.
 \end{equation}
\end{theorem}
\begin{proof} If the operator $\widecheck{T}_{\theta}$ is continuous, Proposition \ref{P8} gives us that $\theta\in (\Lambda_{\infty}(\alpha))^{'}$ and \eqref{Eq6} holds. Let us assume that $\theta\in (\Lambda_{\infty}(\alpha))^{\prime}$ and \eqref{Eq6} holds. Then there exists some $m_{1}\in \mathbb{N}$ and $C_{1}>0$ such that
$$\hspace{1.5in}|\theta_{n-1}|\leq C_{1}e^{m_{1}\alpha_{n}} \hspace{1.55in}\forall n\in \mathbb{N}.$$
From the condition (\ref{Eq6}), for every $k\in \mathbb{N}$, there exist $m_{2}\in \mathbb{N}$ and $C_{2}>0$ such that
\begin{equation*}
\begin{split}
\|\widecheck{T}_{\theta} e_{n}\|_{k}& =\sum^{n}_{j= 1}|\theta_{n-j}|b_{j,k}=\sum^{n-s_{0}}_{j=1} |\theta_{n-j}|b_{j,k} \leq C_{1} \sum^{n-s_{0}}_{j= 1} e^{m_{1}\alpha_{n-j+1}}b_{j,k} \\ & \leq C_{2} \sum^{n-s_{0}}_{j=1} e^{m_{1}\alpha_{n-j+1}}e^{m_{2}\alpha_{j+s_{0}}} \leq C_{2} e^{m_{1}\alpha_{n}} \sum^{n-s_{0}}_{j=1} e^{m_{2}\alpha_{j+s_{0}}}\\ &=C_{2}e^{m_{1}\alpha_{n}}\sum^{n}_{j=1+s_{0}} e^{m_{2}\alpha_{j}}.
\end{split}
\end{equation*}
On the other hand,  $\Lambda_{\infty}(\alpha)$ is nuclear then $\sum^{\infty}_{n=1}e^{-m_{3}\alpha_{n}}$ is convergent for some $m_{3}\in \mathbb{N}$
and then for some $D>0$ we have 
$$ \sum^{n}_{j=1} e^{m_{2}\alpha_{j}}\leq D e^{(m_{2}+m_{3})\alpha_{n}}$$
since 
$$ \sum^{n}_{j=1} e^{m_{2}\alpha_{j}}e^{-(m_{2}+m_{3})\alpha_{n}}\leq \sum^{n}_{j=1} e^{m_{2}\alpha_{j}}e^{-(m_{2}+m_{3})\alpha_{j}}\leq \sum^{\infty}_{j=1} e^{-m_{3}\alpha_{j}} =D<+\infty.$$
Then we can write
\begin{equation*}
\begin{split}
\|\widecheck{T}_{\theta} e_{n}\|_{k} & \leq C_{2} e^{m_{1}\alpha_{n}} \sum^{n}_{j=1+s_{0}} e^{m_{2}\alpha_{j}}\leq  C_{2} e^{m_{1}\alpha_{n}} \sum^{n}_{j=1} e^{m_{2}\alpha_{j}}\\ &\leq D C_{2} e^{m_{1}\alpha_{n}}e^{(m_{2}+m_{3})\alpha_{n}}
=D C_{2}\|e_{n}\|_{m_{1}+m_{2}+m_{3}}.
\end{split}
\end{equation*}
Therefore, $\widecheck{T}_{\theta} e_{n}\in K(b_{n,k})$ for every $n\in \mathbb{N}$ and for all $k\in \mathbb{N}$, there exist $m=m_{1}+m_{2}+m_{3}\in \mathbb{N}$ such that
$$ \sup_{n\in \mathbb{N}}\frac{\|\widecheck{T}_{\theta} e_{n}\|_{k}}{\|e_{n}\|_{m}}< \infty,$$
that is, $\widecheck{T}_{\theta}$ is well defined and continuous by Proposition \ref{P0}.
\end{proof}

We can characterize the compactness of the operators $\widecheck{T}_{\theta} :\Lambda_{\infty}(\beta)\to K(b_{n,k})$ as follows.

\begin{theorem}\label{P10} Let $\Lambda_{\infty}(\alpha)$ be a nuclear power series space of infinite type and $K(b_{n,k})$ be a Montel Köthe space. $\widecheck{T}_{\theta}:\Lambda_{\infty}(\alpha)\to K(b_{n,k})$ is compact if and only if $\theta\in (\Lambda_{\infty}(\alpha))^{\prime}$ and the following condition holds:
\begin{equation}\label{E7}
\hspace{0.2in}\exists m\in \mathbb{N} \quad \forall k\in \mathbb{N} \quad \exists C>0 \quad\quad\quad b_{n-s_{0},k}\leq C e^{m\alpha_{n}}\hspace{0.4in}\forall n>s_{0}.
\end{equation}
\end{theorem}
\begin{proof} If $\widecheck{T}_{\theta}$ is compact, $\theta\in (\Lambda_{\infty}(\alpha))^{\prime}$ from Proposition \ref{P8} and Lemma \ref{Crone} says that there exists a $m\in \mathbb{N}$ such that for all $k\in \mathbb{N}$
$$ \sup_{n}\frac{\|\widecheck{T}_{\theta} e_{n}\|_{k}}{\|e_{n}\|_{m}}< \infty.$$
This gives us that there exists a $C>0$ such that
$$|\theta_{s_{0}}|b_{n-s_{0},k}\leq \|\widecheck{T}_{\theta}e_{n}\|=\sum^{n}_{j=1} |\theta_{n-j}|b_{j,k}\leq C\|e_{n}\|_{m}=Ce^{m\alpha_{n}}.$$
This says that the condition (\ref{E7}) is satisfied.

For the converse, assume  $\theta\in (\Lambda_{\infty}(\alpha))^{\prime}$  and the condition (\ref{E7}) holds. By the same calculation of the proof of the Theorem \ref{P9}, we can write that there exists a $m\in \mathbb{N}$ such that for all $k\in \mathbb{N}$
$$ \sup_{n}\frac{\|\widecheck{T}_{\theta} e_{n}\|_{k}}{\|e_{n}\|_{m}}< \infty.$$
Lemma \ref{Crone} gives us that $\widecheck{T}_{\theta} :\Lambda_{\infty}(\alpha)\to K(b_{n,k})$ is compact.
\end{proof}

\begin{proposition}\label{P11}  Let $\Lambda_{\infty}(\alpha)$ be a nuclear power series space of infinite type and $K(b_{n,k})$ be a Köthe space. 
Assume that there exists a non-decreasing function $S:\mathbb{N}\to \mathbb{N}$ 
such that the following condition holds
\begin{equation}\label{Eq7}
\hspace{0.4in}\forall k\in \mathbb{N} \quad \exists C>0 \quad\quad\quad\quad b_{n,k}\leq C e^{S(k)\alpha_{n}} \hspace{.8in} \forall n\in \mathbb{N}.\end{equation}
Then, the family of the operators $\widecheck{T}_{\theta} :\Lambda_{\infty}(\alpha)\to K(b_{n,k})$ for $\theta\in (\Lambda_{\infty}(\alpha))^{\prime}$, that is,
$$\widecheck{\mathcal{A}}_{\infty}=\{\widecheck{T}_{\theta}:\theta\in (\Lambda_{\infty}(\alpha))^{\prime}\}$$
is 2S-tame.
\end{proposition}
\begin{proof} The proof is similiar to the proof of Theorem \ref{P9}. Let us assume $\theta\in (\Lambda_{\infty}(\alpha))^{\prime}$. Then there exists some $m_{1}\in \mathbb{N}$ and $C_{1}>0$ such that
$$\hspace{1.4in} |\theta_{n-1}|\leq C_{1}e^{m_{1}\alpha_{n}} \hspace{1.4in}\forall n\in \mathbb{N}.$$
Since $\Lambda_{\infty}(\alpha)$ is nuclear then $\sum^{\infty}_{n=1}e^{-m_{2}\alpha_{n}}$ is convergent for some $m_{2}\in \mathbb{N}$ and using similiar to steps in the proof of Theorem \ref{P9} we can show that there exists a $D>0$ such that
$$ \sum^{n}_{j=1} e^{k\alpha_{j}}\leq D e^{(k+m_{2})\alpha_{n}}$$
for every $k\in \mathbb{N}$.
For every $k\in \mathbb{N}$ satisfying $S(k)\geq m_{1}+m_{2}$, there exists $C_{2},C_{3}>0$ such that
\begin{equation*}
\begin{split}
\|\widecheck{T}_{\theta} e_{n}\|_{k}& =\sum^{n}_{j=1}|\theta_{n-j}|b_{j,k} \leq C_{1} \sum^{n}_{j=1} e^{m_{1}\alpha_{n-j+1}}b_{j,k} \\ &\leq C_{2}e^{m_{1}\alpha_{n}} \sum^{n}_{j=1} e^{S(k)\alpha_{j}}\leq C_{3} e^{m_{1}\alpha_{n}} e^{(S(k)+m_{2})\alpha_{n}}  \leq C_{3}e^{2S(k)\alpha_{n}}=\|e_{n}\|_{2S(k)}.
\end{split}
\end{equation*}
Then, we have the following
$$ \|\widecheck{T}_{\theta} e_{n}\|_{k} \leq C_{2}\|e_{n}\|_{2S(k)}$$
for all $k$ satisfying $S(k)\geq m_{1}+m_{2}$, that is, the family $\widecheck{\mathcal{A}}_{\infty}$ is 2S-tame.
\end{proof}
For the converse of Proposition \ref{P10}, we need a stable sequence $\alpha$
when $K(a_{n,k})$ is power series space of finite type $\Lambda_{1}(\alpha)$. 
\begin{theorem}\label{P12} Let $\alpha=\left(\alpha_{n}\right)_{n\in \mathbb{N}}$ be a stable sequence. 
$\widecheck{T}_{\theta} :\Lambda_{1}(\alpha)\to K(b_{n,k})$ is well-defined and continuous if and only if $\theta\in (\Lambda_{1}(\alpha))^{\prime}$ and the following condition holds
\begin{equation}\label{Eq8}
\hspace{0.2in}\forall k\in \mathbb{N} \quad \exists m\in \mathbb{N}, C>0\quad\quad\quad\quad b_{n,k}\leq C e^{-\frac{1}{m}\alpha_{n}} \hspace{0.5in} \forall n\in \mathbb{N}.\end{equation}
\end{theorem}
\begin{proof} Let assume $\theta$ be a sequence and $s_{0}=\min\{t:\theta_{t}\neq 0\}.$ If the operator $\widecheck{T}_{\theta}$ is continuous, Proposition \ref{P8} gives us that $\theta\in (\Lambda_{1}(\alpha))^{'}$ and for every $k\in \mathbb{N}$, there exist $m\in \mathbb{N}$, $C>0$ such that
$$b_{n,k}\leq Ce^{-\frac{1}{m}\alpha_{n+s_{0}}}\leq Ce^{-\frac{1}{m}\alpha_{n}}$$
for every $n\in \mathbb{N}$. This says that 
\eqref{Eq8} holds since $\beta$ is an increasing sequence. For the converse, we assume that assume that $\alpha$ is stable, $\theta\in (\Lambda_{1}(\alpha))^{\prime}$ and \eqref{Eq8} holds. Then there exist some $m_{1}\in \mathbb{N}$ and $C_{1}>0$ such that
$$\hspace{1.45in}|\theta_{n-1}|\leq C_{1}e^{-\frac{1}{m_{1}}\alpha_{n}} \hspace{1.5in} \forall n\in \mathbb{N}.$$
By stability of $\alpha$, we will show that there exists $M>0$ such that for every $n\in \mathbb{N}$ and $j\leq n$ it holds that
\begin{equation}\label{S1} \alpha_{n}\leq M(\alpha_{n-j+1}+\alpha_{j}).
\end{equation}
Since $\alpha$ is stable, there exists a $M>0$ such that $\alpha_{2t}\leq M\alpha_{t}$ for every $t\in \mathbb{N}$. Assume that $n=2t$ or $n=2t+1$ and $1\leq j\leq n$
$$ \alpha_{n}\leq \alpha_{2t+2}\leq M\alpha_{t+1}\leq M (\alpha_{n-j+1}+\alpha_{j})$$
since in this case $t+1\leq j$ or $t+1\leq n-j+1$ and then we have $\alpha_{n-j+1}+\alpha_{j}\geq \alpha_{t+1}$. Therefore the inequality (\ref{S1}) is satisfied.
By using the condition \eqref{Eq8}, we can write that for every $k\in \mathbb{N}$, there exist $m_{2}\in \mathbb{N}$ and $C_{2}>0$ such that
$$ b_{n,k}\leq C_{2}e^{-\frac{1}{m_{2}}\alpha_{n+s_{0}}}\leq C_{2}e^{-\frac{1}{m_{2}}\alpha_{n}}$$
for every $n\in \mathbb{N}$. Let us say $m_{3}=\min(m_{1},m_{2})$. Then we have
\begin{align*}
\begin{split}
\hspace{1in}\|\widecheck{T}_{\theta} e_{n}\|_{k}&=\sum^{n}_{j=1}|\theta_{n-j}|b_{j,k} \leq C_{1}C_{2}\sum^{n}_{j=1}e^{-\frac{1}{m_{1}}\alpha_{n-j+1}}e^{-\frac{1}{m_{2}}\alpha_{j}} \\
&\leq C_{1}C_{2} \sum^{n}_{j=1}e^{-\frac{1}{m_{3}}(\alpha_{n-j+1}+\alpha_{j})} \\ &\leq C_{1}C_{2} \sum^{n}_{j=1}e^{-\frac{1}{Mm_{3}}\alpha_{n}} \leq C_{1}C_{2} n e^{-\frac{1}{m_{3}M}\alpha_{n}}
\end{split}
\end{align*}
Now we choose a $m_{4}\in \mathbb{N}$ so that $m_{4}>m_{3}M$, then 
$$ \hspace{0.5in} -\frac{1}{m_{3}M}+\frac{1}{m_{4}}<0 \hspace{0.5in}\text{and}\hspace{0.5in} \lim_{n\to \infty} ne^{(-\frac{1}{m_{3}M}+\frac{1}{m_{4}})\alpha_{n}}=0.$$
Then there exists a $D>0$ such that for all $n\in \mathbb{N}$
$$ ne^{-\frac{1}{m_{3}M}\alpha_{n}}\leq De^{-\frac{1}{m_{4}}\alpha_{n}}$$
and for some $C_{3}>0$ and for all $n\in \mathbb{N}$
$$ \|\widecheck{T}_{\theta} e_{n}\|_{k} \leq  C_{1}C_{2} n e^{-\frac{1}{m_{3}M}\alpha_{n}}\leq C_{3} 
  e^{-\frac{1}{m_{4}}\alpha_{n}}\leq C_{3} \|e_{n}\|_{m_{4}}.$$
Therefore, for all $k\in \mathbb{N}$, there exist $m\in \mathbb{N}$ such that
$$ \sup_{n\in \mathbb{N}}\frac{\|\widecheck{T}_{\theta}e_{n}\|_{k}}{\|e_{n}\|_{m}}< \infty,$$
that is, $\widecheck{T}_{\theta} :\Lambda_{1}(\alpha)\to K(b_{n,k})$ is well-defined and  continuous  by Proposition \ref{P0}.
\end{proof}

The following theorem characterize the compactness of the operators $\widecheck{T}_{\theta}:\Lambda_{1}(\alpha)\to K(b_{n,k})$.

\begin{theorem}\label{P13} Let $\alpha=(\alpha_{n})_{n\in \mathbb{N}}$ be a stable sequence and $K(b_{n,k})$ be a Montel Köthe space. $\widecheck{T}_{\theta}:\Lambda_{1}(\alpha)\to K(b_{n,k})$ is compact if and only if $\theta\in (\Lambda_{1}(\alpha))^{\prime}$ and the following condition holds:
\begin{equation}\label{Eq9}
\hspace{.2in} \exists m\in \mathbb{N} \quad \forall k\in \mathbb{N} \quad \exists C>0 \quad\quad\quad b_{n,k}\leq C e^{-\frac{1}{m}\alpha_{n}} \hspace{.5in}\forall n\geq k.
\end{equation}
\end{theorem}
\begin{proof} Let assume $\theta$ be a sequence and $s_{0}=\min\{t:\theta_{t}\neq 0\}.$ If the operator $\widehat{T}_{\theta}$ is continuous, then $\theta\in (\Lambda_{1}(\alpha))^{\prime}$ from Proposition \ref{P8} and Proposition \ref{P8C} says that for every $k\in \mathbb{N}$, there exist $m\in \mathbb{N}$, $C>0$ such that
$$b_{n,k}\leq Ce^{-\frac{1}{m}\alpha_{n+s_{0}}} \leq Ce^{-\frac{1}{m} \alpha_{n}}$$
for every $n\in \mathbb{N}$. This says that 
\eqref{Eq9} holds since $\alpha$ is an increasing sequence.

For the converse, assume that  $\theta\in (\Lambda_{1}(\alpha))^{\prime}$ and the condition (\ref{Eq9}) holds. By the same calculation of the proof of the Theorem \ref{P12}, we can write that there exists  $m\in \mathbb{N}$  such that
for all $k\in \mathbb{N}$
$$ \sup_{n}\frac{\|\widecheck{T}_{\theta}e_{n}\|_{k}}{\|e_{n}\|_{m}}< \infty,$$
that is, $\widecheck{T}_{\theta}:\Lambda_{1}(\alpha) \to K(a_{n,k})$ is compact from Lemma \ref{Crone}.
\end{proof}
\begin{proposition}\label{P14} Let $\alpha=(\alpha_{n})_{n\in \mathbb{N}}$ be a stable sequence and $\displaystyle M=\sup_{n\in\mathbb{N}}\frac{\alpha_{2n}}{\alpha_{n}}$. Assume that there exists a non-decreasing function $S:\mathbb{N}\to \mathbb{N}$ such that the following condition hold:
\begin{equation}\label{E10}
\hspace{0.3in}\forall k\in \mathbb{N} \quad \exists C>0 \quad\quad \quad\quad\quad b_{n,k}\leq C e^{-\frac{1}{S(k)}\alpha_{n}} \hspace{0.6in}\forall n\in \mathbb{N}.
\end{equation}
Then, the family of the operators $\widecheck{T}_{\theta} :\Lambda_{1}(\alpha)\to K(b_{n,k})$ for $\theta\in (\Lambda_{1}(\beta))^{\prime}$, that is,
$$\widecheck{\mathcal{A}}_{1}=\{\widecheck{T}_{\theta}:\theta\in (\Lambda_{1}(\alpha))^{\prime}\}$$
is $\tilde{S}$-tame where $\tilde{S}(k)=(M+1)S(k)$ for every $k\in \mathbb{N}$.
\end{proposition}
\begin{proof} The proof is similar to the proof of the Theorem \ref{P12}. Since $\theta\in (\Lambda_{1}(\alpha))^{\prime}$, there exists some $m_{1}\in \mathbb{N}$ and $C_{1}>0$ such that
$$\hspace{1.4in}|\theta_{n-1}|\leq C_{1}e^{-\frac{1}{m_{1}}\alpha_{n}} \hspace{1.4in} \forall n\in \mathbb{N}.$$
By using \eqref{E10}, we can write that for every $k\in \mathbb{N}$ satisfying $S(k)\geq m_{1}$, there exists a $C_{2},C_{3}, C_{4}>0$ such that
\begin{align*}
\begin{split}
\|\widecheck{T}_{\theta}e_{n}\|_{k}&=\sum^{n}_{j=1}|\theta_{n-j}|b_{j,k} \leq C_{1} \sum^{n}_{j=1} e^{-\frac{1}{m_{1}}\alpha_{n-j+1}}b_{j,k}
\leq C_{2} \sum^{n}_{j=1} e^{-\frac{1}{m_{1}}\alpha_{n-j+1}}e^{-\frac{1}{S(k)}\alpha_{j}} \\ &\leq C_{2} \sum^{n}_{j=1}e^{-\frac{1}{S(k)}(\alpha_{n-j+1}+\alpha_{n})} \leq C_{3}\sum^{n}_{j=1} e^{-\frac{1}{MS(k)}\alpha_{n}}=C_{3}n e^{-\frac{1}{MS(k)}\alpha_{n}}  \\ &  \leq C_{4}e^{-\frac{1}{(M+1)S(k)}\alpha_{n}} \leq  C_{4}\|e_{n}\|_{(M+1)S(k)}
\end{split}
\end{align*}
Then for all $k\in \mathbb{N}$ satisfying $S(k)\geq m_{1}$, then, we have the following 
$$\|\widecheck{T}_{\theta} e_{n}\|_{k}\leq C_{3}\|e_{n}\|_{(M+1)S(k)}$$
that is the family $\widecheck{\mathcal{A}}_{1}$ is $\tilde{S}$-tame, where $\tilde{S}(k)=(M+1)S(k)$ for all $k\in \mathbb{N}$.
\end{proof}
\section{Toeplitz Operators Defined between Power Series Spaces}
When given a Toeplitz matrix
$$
\begin{pmatrix}
\theta_{0} &\theta_{-1}&\theta_{-2}&\cdots \\
\theta_{1}&\theta_{0}&\theta_{-1}&\cdots\\ \theta_{2}&\theta_{1}&\theta_{0}&\cdots\\
\vdots&\vdots&\vdots&\ddots
\end{pmatrix},
$$
one can express this matrix as a sum of a lower and an upper Toeplitz matrix in the following way:
$$
\begin{pmatrix}
\theta_{0} &\theta_{-1}&\theta_{-2}&\cdots \\
\theta_{1}&\theta_{0}&\theta_{-1}&\cdots\\ \theta_{2}&\theta_{1}&\theta_{0}&\cdots\\
\vdots&\vdots&\vdots&\ddots
\end{pmatrix} =
\begin{pmatrix}
\theta^{\prime}_{0} &0&0&\cdots \\
\theta_{1}&\theta^{\prime}_{0}&0&\cdots\\ \theta_{2}&\theta_{1}&\theta^{\prime}_{0} &\cdots\\
\vdots&\vdots&\vdots&\ddots
\end{pmatrix}
+
\begin{pmatrix}
 \theta^{\prime\prime}_{0}&\theta_{-1}&\theta_{-2}&\cdots \\
0&\theta^{\prime\prime}_{0}&\theta_{-1}&\cdots\\ 0&0&\theta^{\prime\prime}_{0}&\cdots\\
\vdots&\vdots&\vdots&\ddots
\end{pmatrix}
$$
Here $\theta_{0}=\theta^{\prime}_{0}+ \theta^{\prime\prime}_{0}$ and we can choose $\theta^{\prime}_{0}$ and $\theta^{\prime\prime}_{0}$  to be non-zero.
Let define the sequences
$$\widehat{\theta}=(\theta^{\prime}_{0},\theta_{1},\theta_{2},\theta_{3},\dots)$$
and 
$$ \widecheck{\theta}= (\theta^{\prime\prime}_{0},\theta_{-1},\theta_{-2},\theta_{-3},\dots).$$
Therefore every Toeplitz matrix can be associated with two sequences. Denoting the operator defined by the Toeplitz matrix as $T_{\theta}$, one can write this operator as a sum of two operators defined by a lower and an upper Toeplitz matrix, that is,
$$T_{\theta}=\widehat{T}_{\widehat{\theta}}+\widecheck{T}_{\widecheck{\theta}}.$$

In this section, we collect some results regarding the continuity and compactness of a Toeplitz operator by utilizing the theorems proven in the previous sections. 
\begin{theorem}\label{T11} Let $\alpha$ be a stable sequence. A Toeplitz operator $T_{\theta}: \Lambda_{1}(\alpha)\to \Lambda_{1}(\beta)$  is well-defined and  continuous if  $\theta\in \Lambda_{1}(\beta)\oplus (\Lambda_{1}(\alpha))^{\prime}$ and the following condition holds:
\begin{equation}\label{Et1}
\hspace{.2in} \forall k\in \mathbb{N} \quad \exists m\in \mathbb{N}, C>0 \quad\quad\quad e^{-\frac{1}{k}\beta_{n}}\leq C e^{-\frac{1}{m}\alpha_{n}} \hspace{.6in} \forall n\in \mathbb{N}.\end{equation}
\end{theorem}
\begin{proof} Let assume that Toeplitz matrix associated with $T_{\theta}$ is given by sequnces $\widehat{\theta}$ and $\widecheck{\theta}$. We can choose $\theta^{\prime}_{0}$ and $\theta^{\prime\prime}_{0}$  to be non-zero. Then Toeplitz operator $T_{\theta}$ is defined as
$$ T_{\theta}=\widehat{T}_{\widehat{\theta}}+\widecheck{T}_{\widecheck{\theta}}.$$
From Theorem \eqref{P2}, $\widehat{T}_{\widehat{\theta}}$ is well-defined and continous if and only if $\widehat{\theta}\in \Lambda_{1}(\beta)$ and the condition \eqref{Et1} is satisfied. From Theorem \ref{P12}, $\widecheck{T}_{\widecheck{\theta}}$ is well-defined and  continous if and only if $\widecheck{\theta}\in (\Lambda_{1}(\alpha))^{\prime}$ and the condition \eqref{Et1} is satisfied. Therefore, these give us that $T_{\theta}$ is well-defined  and continuous if $\theta=\widehat{\theta}\oplus \widecheck{\theta}\in \Lambda_{1}(\beta)\oplus (\Lambda_{1}(\alpha))^{\prime}$ and the condition \eqref{Et1} is satisfied.
\end{proof}

\begin{theorem} Let $\alpha$ be a stable sequence. A Toeplitz operator $T_{\theta}: \Lambda_{1}(\alpha)\to \Lambda_{1}(\beta)$  is compact if a $\theta\in \Lambda_{1}(\beta)\oplus (\Lambda_{1}(\alpha))^{\prime}$ and the following condition holds:
\begin{equation}\label{E11}
\hspace{.2in} \exists m\in \mathbb{N} \quad \forall k\in \mathbb{N} \quad \exists C>0 \quad\quad\quad e^{-\frac{1}{k}\beta_{n}}\leq C e^{-\frac{1}{m}\alpha_{n}} \hspace{.4in} \forall n\geq k.
\end{equation}
 \end{theorem}
\begin{proof} The proof is similar to the proof of Theorem \ref{T11} and follows from Theorem \ref{P3} and Theorem \ref{P13}.
\end{proof}
\begin{rem} The condition in Theorem \ref{E11} is satisfied for the sequences $\alpha=(n)_{n\in \mathbb{N}}$ and $\beta=(n^{2})_{n\in \mathbb{N}}$.
\end{rem}

\begin{theorem} 
Let $\beta$ be a stable sequence and $\Lambda_{\infty}(\alpha)$ be a nuclear power series space of infinite type.  A Toeplitz operator $T_{\theta}: \Lambda_{\infty}(\alpha)\to \Lambda_{\infty}(\beta)$ is well-defined and  continuous if $\theta\in \Lambda_{\infty}(\beta)\oplus (\Lambda_{\infty}(\alpha))^{\prime}$ and the following condition holds:
\begin{equation}\label{}
\hspace{.2in}\forall k\in \mathbb{N} \quad \exists m\in \mathbb{N}, C>0 \quad\quad\quad\quad e^{k\beta_{n}}\leq C e^{m\alpha_{n}}\hspace{.3in} \forall n\geq k.\end{equation}
\end{theorem}
\begin{proof} This follows from Theorem \ref{P5} and Theorem \ref{P9}.
\end{proof}
\begin{theorem} \label{P15}
Let $\beta$ be a stable sequence and $\Lambda_{\infty}(\alpha)$ be a nuclear power series space of infinite type. A Toeplitz operator $T_{\theta}: \Lambda_{\infty}(\alpha)\to \Lambda_{\infty}(\beta)$ is compact if $\theta\in \Lambda_{\infty}(\beta)\oplus (\Lambda_{\infty}(\alpha))^{\prime}$ and the following condition holds:
\begin{equation}\label{E12}
\hspace{.2in}\exists m\in \mathbb{N} \quad \forall k\in \mathbb{N} \quad \exists C>0 \quad\quad\quad\quad e^{k\beta_{n}}\leq C e^{m\alpha_{n}} \hspace{.3in}\forall n\geq k.
\end{equation}
\end{theorem}
\begin{proof} This follows from Theorem \ref{P6} and Theorem \ref{P10}.
\end{proof}
\begin{rem} The condition in Theorem \ref{P15} is satisfied for the sequences $\alpha=(n^{2})_{n\in \mathbb{N}}$ and $\beta=(n)_{n\in \mathbb{N}}$.
\end{rem}



\begin{theorem} Let 
$\Lambda_{\infty}(\alpha)$ be a nuclear power series space of infinite type and $\Lambda_{1}(\beta)$ be a nuclear power series spaces of finite type.  $T_{\theta}: \Lambda_{\infty}(\alpha)\to \Lambda_{1}(\beta)$ is  well-defined, continuous and compact if $\theta\in \Lambda_{1}(\beta)\oplus (\Lambda_{\infty}(\beta))^{\prime}$.
\end{theorem}
\begin{proof} Continuity follows from Theorem \ref{P2} and Theorem \ref{P9}. Compactness follows from Theorem \ref{P3} and Theorem \ref{P10}.
\end{proof}


\begin{rem} By considering the conditions in Proposition \ref{P1} or  Proposition \ref{P8}, the Toeplitz operator $T_{\theta}$ is not well defined from $\Lambda_{1}(\alpha)$ to $\Lambda_{\infty}(\beta)$. 
\end{rem}
 
\begin{center}
Acknowledgments
\end{center}
The results in this paper were obtained while the author visited at University of Toledo. I would like to thank TUBITAK for their support and Prof. Dr. Sönmez Şahutoğlu for sharing Jasiczak's paper with me, which enabled me to create this work.

\end{document}